\documentclass[11pt]{article}
\UseRawInputEncoding
\usepackage{mathrsfs}
\usepackage{amssymb}
\usepackage{amsmath}
\usepackage{amsbsy}
\usepackage{epsfig}
\usepackage{enumerate}
\usepackage{bm}
\usepackage{color}
\usepackage[english]{babel}

\usepackage{bbm}
\usepackage{mathrsfs}
\usepackage{dsfont}
\usepackage{graphicx,theorem}

\makeatletter

\newcommand{\Rmnum}[1]{\expandafter\@slowromancap\romannumeral #1@}
\makeatother

 \newcommand{\me}{\mathrm{e}} 
\newcommand{\dif}{\mathrm{d}} \DeclareMathAlphabet{\mathsfsl}{OT1}{cmss}{m}{sl} \DeclareMathAlphabet{\mathpzc}{OT1}{pzc}{m}{it}

    \newcommand{\ee}{\mathbb{E}}
\newcommand{\hh}{\mathbb{H}}   \newcommand{\nn}{\mathbb{N}} \newcommand{\rr}{\mathbb{R}}
    \newcommand{\vv}{\mathbb{V}}

    \def\FF{\mathcal F}  \def\HH{\mathcal H}

 \def\d"{^{\prime\prime}} \def\d'{^{\prime}}

\theoremstyle{plain}
\newtheorem{lem}{Lemma}[section]  \newtheorem{thm}{Theorem}[section]
\newtheorem{cor}{Corollary}[section] \newtheorem{defn}{Definition}[section] \newtheorem{rmk}{Remark}[section] 
\numberwithin{equation}{section}

\renewcommand{\baselinestretch}{1.5}

\begin{document}

\begin{center}{\LARGE\bf Equivalent conditions of complete convergence for weighted sums of sequences of i. i. d. random variables under sublinear expectations}
\end{center}

\smallskip
%\centerline{\textit{April 2021}}
\begin{center} {\sc
Mingzhou Xu\footnote{This work was supported by Supported by Doctoral Scientific Research Starting Foundation of Jingdezhen Ceramic University ( Nos.102/01003002031 ), Scientific Program of Department of Education of Jiangxi Province of China (Nos. GJJ190732, GJJ180737), Natural Science Foundation Program of Jiangxi Province 20202BABL211005, and National Natural Science Foundation of China (Nos. 61662037)
} \quad Kun Cheng
}\\
{\sl \small School of Information Engineering, Jingdezhen Ceramic University
  Jingdezhen 333403, China} \\
(Email:mingzhouxu@whu.edu.cn\quad chengkun0010@126.com) \\

\end{center}

 \renewcommand{\abstractname}{~}
\begin{abstract}
{\bf Abstract:}
The complete convergence for weighted sums of sequences of independent, identically distributed random variables under sublinear expectations space was studied. By moment inequality and truncation methods, we establish the equivalent conditions of complete convergence for weighted sums of sequences of independent, identically distributed random variables under sublinear expectations space. The results extend the corresponding results obtained by Guo (2012) to those for sequences of independent, identically distributed random variables under sublinear expectations space.



% {\bf Short text for running head:}  LLN  under sub-linear expectations
\vspace{-3mm}
\end{abstract}

%%%%%%%%%%%%%%%%%%%%%%%%%%%%%%%%%%%%%%%%%%%%%%%%%%%%%%%%%%%%
\renewcommand{\baselinestretch}{1.2}

%%%%%%%%%%%%%%%%%%%%%%%%%%%%%%%%%%%%%%%%%%%%%%%%%%%%%%%%%%%%
%% Text of article.
%%%%%%%%%%%%%%%%%%%%%%%%%%%%%%%%%%%%%%%%%%%%%%%%%%%%%%%%%%%%
%    Section headings
%\baselineskip 11pt\parindent=10.8pt

\section{ Introduction and notations.}\label{sect1}
\setcounter{equation}{0}

Motivated by the study of the uncertainty of probability and distribution, Peng \cite{Peng2007,Peng2010} introduced the concept of  the sublinear expectation space. We consider the sublinear expectation as an extension of linear expectation. Inspired by the work of Peng \cite{Peng2007,Peng2010}, people try to study lots of limit theorems under sublinear expectation space to extend the corresponding results in probability and statistics. Zhang \cite{Zhang2015,Zhang2016a,Zhang2016b} obained the exponential inequalities, Rosenthal's inequalities, Donsker's invariance principle under sublinear expectations. Under sublinear expectations, Wu \cite{Wu2020} established precise asymptotics for complete integral convergence, Xu and Zhang \cite{Xujiapan2019,Xujiapan2020} proved the law of logarithm for arrays of random variables and three series theorem for independent random variables, Chen \cite{Chen2016} obtained strong laws of large numbers, Xu and Cheng\cite{Xu2021a,Xu2021c} investigated precise asymptotics in the law of the iterated logarithm and convergence for sums of  independent, identically distributed random variables. Under sublinear expectations  the interested reader could refer to Gao and Xu \cite{Gao2011}, Fang et al. \cite{Fang12018},  Hu et al. \cite{Hufeng2014}, Kuczmaszewska \cite{Kuczm}, Wang and Wu \cite{Wangwenjuan2019}, Hu and Yang \cite{Huzechun2017}, Zhang \cite{Zhang2016c}, Yu and Wu \cite{Yudonglin2018}, Wu and Jiang \cite{Wuqunying2018}, Ma and Wu \cite{Maxiaochen2020}, Huang and Wu \cite{Huang2019}, Zhong and Wu \cite{Zhong2017} and references therein for more limit theorems.

Meng et al. \cite{Meng2019} studied convergence for sums of asymptotically almost negatively associated random variables. Guo \cite{Guo2012} obtained  equivalent conditions of complete convergence for weighted sums of sequences of negatively dependent random variables. Guo and Shan \cite{Guo2020} established equivalent conditions of complete $q$th moment convergence for weighted sums of sequences of negatively orthant dependent random variables.   For references on complete convergence in linear expectation space, the reader could refer to Meng et al. \cite{Meng2019}, Guo \cite{Guo2012}, Shen and Wu \cite{Shen2014} and references therein.  In this paper, we investigate the complete convergence for independent, identically distributed (i. i. d.) random variables under sublinear expectations space. We proved equivalent conditions of complete convergence for weighted sums of sequences of independent, identically distributed random variables under sublinear expectations space.

The rest of this paper is organized as follows. In the next section, we summarize necessary basic notions, concepts and relevant properties, and give necessary lemmas under sublinear expectations. In Section 3, we present our main results, Theorems \ref{thm1}, \ref{thm2},  the proofs of which are given in Section 4.

%%%%%%%%%%%%%%%%%%%%%%%%%%%%%%%%%%%%%%%%%%%%%%%%%%%

\section{Preliminaries}
\setcounter{equation}{0}
We use notations similar to that of Peng \cite{Peng2010}. Assume that $(\Omega,\FF)$ is a given measurable space. Suppose that $\HH$ is a subset of all random variables on $(\Omega,\FF)$ such that $I_A\in \HH$, where $A\in\FF$ and  $X_1,\cdots,X_n\in \HH$ implies $\varphi(X_1,\cdots,X_n)\in \HH$ for each $\varphi\in C_{l,Lip}(\rr^n)$, where $C_{l,Lip}(\rr^n)$ denotes the linear space of (local lipschitz) function $\varphi$ satisfying
$$
|\varphi(\mathbf{x})-\varphi(\mathbf{y})|\le C(1+|\mathbf{x}|^m+|\mathbf{y}|^m)(|\mathbf{x}-\mathbf{y}|), \forall \mathbf{x},\mathbf{y}\in \rr^n
$$
for some $C>0$, $m\in \nn$ depending on $\varphi$. We regard $\HH$  as the space of random variables.
\begin{defn}\label{defn1}({\em See Peng \cite{Peng2007,Peng2010}}) A sublinear expectation $\ee$ on $\HH$ is a functional $\ee:\hh\mapsto \bar{\rr}:=[-\infty,\infty]$ satisfying the following properties: for all $X,Y\in \HH$, we have
\begin{description}
  \item[\rm (a)]  Monotonicity: If $X\ge Y$, then $\ee[X]\ge \ee[Y]$;
\item[\rm (b)] Constant preserving: $\ee[c]=c$, $\forall c\in\rr$;
\item[\rm (c)] Positive homogeneity: $\ee[\lambda X]=\lambda\ee[X]$, $\forall \lambda\ge 0$;
\item[\rm (d)] Sub-additivity: $\ee[X+Y]\le \ee[X]+\ee[Y]$ whenever $\ee[X]+\ee[Y]$ is not of the form $\infty-\infty$ or $-\infty+\infty$.
 \end{description}
\end{defn}

A set function $V:\FF\mapsto[0,1]$ is named a capacity if it satisfies
\begin{description}
\item[\rm (a)]$V(\emptyset)=0$, $V(\Omega)=1$;
\item[\rm (b)]$V(A)\le V(B)$, $A\subset B$, $A,B\in \FF$.
 \end{description}
A capacity $V$ is said to be sub-additive if it satisfies $V(A+B)\le V(A)+V(B)$, $A,B\in \FF$.

In this paper, given a sublinear expectation space $(\Omega, \HH, \ee)$, we define a capacity: $\vv(A):=\inf\{\ee[\xi]:I_A\le \xi, \xi\in \HH\}$, $\forall A\in \FF$ (see Zhang \cite{Zhang2016a}). Obviously $\vv$ is a sub-additive capacity. We define the Choquet expectations $C_{\vv}$ by
$$
C_{\vv}(X):=\int_{0}^{\infty}\vv(X>x)\dif x +\int_{-\infty}^{0}(\vv(X>x)-1)\dif x.
$$

Assume that $\mathbf{X}=(X_1,\cdots, X_m)$, $X_i\in\HH$ and $\mathbf{Y}=(Y_1,\cdots,Y_n)$, $Y_i\in \HH$  are two random variables on  $(\Omega_1, \HH, \ee)$. $\mathbf{Y}$ is said to be independent of $\mathbf{X}$, if for each Borel-measurable function $\psi$ on $\rr^m\times \rr^n$ with $\psi(\mathbf{X},\mathbf{Y}), \psi(\mathbf{x},\mathbf{Y})\in \HH$ for each $x\in\rr^m$, we have $\ee[\psi(\mathbf{X},\mathbf{Y})]=\ee[\ee\psi(\mathbf{x},\mathbf{Y})|_{\mathbf{x}=\mathbf{X}}]$ whenever $\bar{\psi}(\mathbf{x}):=\ee[|\psi(\mathbf{x},\mathbf{Y})|]<\infty$ for each $\mathbf{x}$ and $\ee[|\bar{\psi}(\mathbf{X})|]<\infty$ (cf. Definition 2.5 in Chen \cite{Chen2016} ). $\{X_n\}_{n=1}^{\infty}$ is said to be a sequence of independent random variables, if $X_{n+1}$ is independent of $(X_1,\cdots,X_n)$ for each $n\ge 1$.

Suppose that $\mathbf{X}_1$ and $\mathbf{X}_2$ are two $n$-dimensional random vectors defined, respectively, in sublinear expectation spaces $(\Omega_1,\HH_1,\ee_1)$ and $(\Omega_2,\HH_2,\ee_2)$. They are said to be identically distributed if  for every Borel-measurable function $\psi$ such that $\psi(X_1), \psi(X_2)\in \HH$,
$$
\ee_1[\psi(\mathbf{X}_1)]=\ee_2[\psi(\mathbf{X}_2)], \mbox{  }
$$
whenever the sublinear expectations are finite. $\{X_n\}_{n=1}^{\infty}$ is said to be identically distributed if for each $i\ge 1$, $X_i$ and $X_1$ are identically distributed.

Throughout this paper we assume that $\ee$ is countably sub-additive, i.e., $\ee(X)\le \sum_{n=1}^{\infty}\ee(X_n)$, whenever $X\le \sum_{n=1}^{\infty}X_n$, $X,X_n\in \HH$, and $X\ge 0$, $X_n\ge 0$, $n=1,2,\ldots$. We denote by $C$ a positive constant which may vary from line to line. $I(A)$ denotes the indicator function of $A$, $a_n\ll b_n$ means that there exists a constant $C>0$ such that $a_n\le c b_n$ for $n$ large enough and let $a_n\approx b_n$ denote that $a_n\ll b_n$ and $b_n\ll a_n$. Let $\log x$ denote $\ln\max\{\me x\}$.

To establish our results, we present the following lemmas.
\begin{lem}\label{lem1} ( See Lemma 2.1 in Xu and Cheng \cite{Xu2021c})Let $X$ be a random variable under sublinear expectation space $(\Omega,\HH,\ee)$ . Then for any $\alpha> 0$, $\gamma>0$, and $\beta>-1$,
\begin{eqnarray*}
&(\rm\expandafter{\romannumeral1})&\int_{1}^{\infty}u^{\beta}C_{\vv}(|X|^{\alpha}I(|X|>u^{\gamma}))\dif u\le C C_{\vv}(|X|^{(\beta+1)/\gamma+\alpha}),\\
&(\rm\expandafter{\romannumeral2})&\int_{1}^{\infty}u^{\beta}\ln(u)C_{\vv}(|X|^{\alpha}I(|X|>u^{\gamma}))\dif u\le C C_{\vv}(|X|^{(\beta+1)/\gamma+\alpha}\ln(1+|X|)).
 \end{eqnarray*}
\end{lem}
Write $S_k=X_1+\cdots+X_k$, $S_0=0$.
\begin{lem}\label{lem2} (cf. Corollary 2.2 and Theorem 2.3 in Zhang \cite{Zhang2016b}) Suppose that $X_{k+1}$ is independent of $(X_1,\ldots,X_k)$ under sublinear expectation space $(\Omega,\HH,\ee)$ with $\ee(X_i)\le 0$, $k=1,\ldots,n-1$. Then
\begin{equation}\label{01}
\ee\left[\left|\max_{k\le n}(S_n-S_k)\right|^p\right]\le 2^{2-p}\sum_{k=1}^{n}\ee[|X_k|^p], \mbox{   for $1\le p\le 2$,}
\end{equation}
\begin{equation}\label{02}
\ee\left[\left|\max_{k\le n}(S_n-S_k)\right|^p\right]\le C_p\left\{\sum_{k=1}^{n}\ee[|X_k|^p]+\left(\sum_{k=1}^{n}\ee[|X_k|^2\right)^{p/2}\right\}, \mbox{   for $ p\ge 2$.}
\end{equation}
\end{lem}
By Lemma \ref{lem2} and the proofs of Theorem 3 in M\'{o}ricz \cite{Moricz1976} combined with Minkowski's inequality under sublinear expectations, we can get the following.
\begin{lem}\label{lem3} Let $\{X_n;n\ge 1\}$ be a sequence of independent random variables under sublinear expectation space $(\Omega,\HH,\ee)$ with $\ee(X_i)\le 0$, $k=1,\ldots,n-1$.. Then for $M\ge 2$,
\begin{eqnarray}\label{03}
\ee\max_{1\le j\le n}\left|\sum_{i=1}^{j}X_i\right|^M\le C\log^M n\left(\sum_{i=1}^{n}\ee|X_i|^M+\left(\sum_{i=1}^{n}\ee|X_i|^2\right)^{M/2}\right).
\end{eqnarray}
\end{lem}

\begin{lem}\label{lem4}(See Lemma 2.4 in Xu and Cheng \cite{Xu2021c}) Let $\{X_n;n\ge 1\}$ be a sequence of independent random variables under sublinear expectation space $(\Omega,\HH,\ee)$. Then for all $n\ge 1$ and $x>0$,
\begin{eqnarray}\label{04}
&&\left[1-\vv\left(\max_{1\le j\le n}|X_j|>x\right)\right]^2\sum_{j=1}^{n}\vv(|X_j|>x)\le 4\vv\left(\max_{1\le j\le n}|X_j|>x\right).
\end{eqnarray}
\end{lem}

\section{Main results}
\setcounter{equation}{0}
Our main results are as follows.
\begin{thm}\label{thm1} Let $\{X_n,n\ge 1\}$ be a sequence of independent random variables, identically distributed as $X$ under sublinear expectation space  $(\Omega,\HH,\ee)$. Assume that $r>1$, $p>\frac12$, $\beta+p>0$, and suppose that $\ee X=-\ee(-X)=0$ for $\frac12<p\le 1$. Suppose that $\{a_{ni}\approx (i/n)^{\beta}(1/n)^p, 1\le i\le n, n\ge 1\}$ is a triangular array of real numbers. Then the following is equivalent:
\begin{description}
 \item[\rm (i)]\begin{eqnarray}\label{06}
\begin{cases} C_{\vv}\left(|X|^{(r+1)/(p+\beta)}\right)<\infty,& \text{
for $-p<\beta<-p/r$,}\\
C_{\vv}\left(|X|^{r/p}\ln(1+|X|)\right)<\infty, & \text{ for $\beta=-p/r$,}\\
C_{\vv}\left(|X|^{r/p}\right)<\infty, & \text{ for $\beta>-p/r$,}
\end{cases}
\end{eqnarray}
\item[\rm (ii)]
\begin{equation}\label{07}
\sum_{n=1}^{\infty}n^{r-2}\vv\left(\max_{1\le j\le n}\left|\sum_{i=1}^{j}a_{ni}X_{i}\right|>\epsilon\right)<\infty, \mbox{ $\forall \epsilon>0$}
\end{equation}
\end{description}
\end{thm}
\begin{thm}\label{thm2} Let $\{X_n,n\ge 1\}$ be a sequence of independent random variables, identically distributed as $X$ under sublinear expectation space  $(\Omega,\HH,\ee)$. Assume that $r>1$, $p>\frac12$, $\beta+p>0$, and suppose that $\ee X=-\ee(-X)=0$ for $\frac12<p\le 1$. Suppose that $\{a_{ni}\approx ((n-i)/n)^{\beta}(1/n)^p, 0\le i\le n-1, n\ge 1\}$ is a triangular array of real numbers. Then the following are equivalent:
\begin{description}
 \item[\rm (i)]\begin{eqnarray}\label{08}
\begin{cases} C_{\vv}\left(|X|^{(r-1)/(p+\beta)}\right)<\infty,& \text{
for $-p<\beta<-p/r$,}\\
C_{\vv}\left(|X|^{r/p}\ln(1+|X|)\right)<\infty, & \text{ for $\beta=-p/r$,}\\
C_{\vv}\left(|X|^{r/p}\right)<\infty, & \text{ for $\beta>-p/r$,}
\end{cases}
\end{eqnarray}
\item[\rm (ii)]
\begin{equation}\label{09}
\sum_{n=1}^{\infty}n^{r-2}\vv\left(\max_{0\le j\le n-1}\left|\sum_{i=1}^{j}a_{ni}X_{i}\right|>\epsilon\right)<\infty, \mbox{ $\forall \epsilon>0$}
\end{equation}
\end{description}
\end{thm}
\begin{rmk} Take $p=1$ in Theorem \ref{thm1}, we obtain Theorem 3.1 in Xu and Cheng \cite{Xu2021b,Xu2021c}.
\end{rmk}

\begin{cor}\label{cor1} Let $\{X_n,n\ge 0\}$ be a sequence of independent random variables, identically distributed as $X$  under sublinear expectation space  $(\Omega,\HH,\ee)$ with $\ee(X_i)=-\ee(-X_i)= 0$, $i=1,2,\ldots$. Assume that $r>1$, $p>\frac12$, $0<\alpha\le 1$. Let $A_n^{\alpha}=(\alpha+1)(\alpha+2)\cdots(\alpha+n)/n!$, $n=1,2,\ldots$ and $A_0^{\alpha}=0$.  The following are equivalent:
\begin{description}
 \item[\rm (i)]\begin{eqnarray}\label{10}
\begin{cases} C_{\vv}\left(|X|^{(r-1)/(p\alpha)}\right)<\infty,& \text{
for $0<\alpha<1-1/r$,}\\
C_{\vv}\left(|X|^{r/p}\ln(1+|X|)\right)<\infty, & \text{ for $\alpha=1-1/r$,}\\
C_{\vv}\left(|X|^{r/p}\right)<\infty, & \text{ for $1-1/r<\alpha\le 1$,}
\end{cases}
\end{eqnarray}
\item[\rm (ii)]
\begin{equation}\label{11}
\sum_{n=1}^{\infty}n^{r-2}\vv\left(\max_{0\le j\le n-1}\left|\sum_{i=0}^{j}(A_{n-i}^{\alpha-1})^pX_{i}\right|>\epsilon(A_n^{\alpha})^p\right)<\infty, \mbox{ $\forall \epsilon>0$}
\end{equation}
\end{description}
\end{cor}

\section{Proofs of Theorems \ref{thm1},\ref{thm2}}

{\bf Proof}. [Proof of Theorem \ref{thm1}]  First we prove (\ref{06}) $\Rightarrow$ (\ref{07}). Choose $\delta>0$ small enough and sufficiently large integer $K$. For all $1\le i\le n$, $n\ge 1$, write
\begin{eqnarray}\label{12}
\nonumber &&X_{ni}^{(1)}=-n^{-\delta}I(a_{ni}X_i<-n^{-\delta})+a_{ni}X_iI(|a_{ni}X_i|\le n^{-\delta})+n^{-\delta}I(a_{ni}X_i>n^{-\delta}),\\
\nonumber &&X_{ni}^{(2)}=(a_{ni}X_i-n^{-\delta})I(n^{-\delta}<a_{ni}X_i<\frac{\epsilon}{K}),\\
&&X_{ni}^{(3)}=(a_{ni}X_i+n^{-\delta})I(-\frac{\epsilon}{K}<a_{ni}X_i<-n^{-\delta}),\\
\nonumber &&X_{ni}^{(3)}=(a_{ni}X_i+n^{-\delta})I(a_{ni}X_i\le -\frac{\epsilon}{K})+(a_{ni}X_i-n^{-\delta})I(a_{ni}X_i\ge\frac{\epsilon}{K}).
 \end{eqnarray}
 Obviously, $\sum_{i=1}^{k}a_{ni}X_i=\sum_{i=1}^{k}X_{ni}^{(1)}+\sum_{i=1}^{k}X_{ni}^{(2)}+\sum_{i=1}^{k}X_{ni}^{(3)}+\sum_{i=1}^{k}X_{ni}^{(4)}$. Notice that
 \begin{eqnarray}\label{13}
\left(\max_{1\le k\le n}\left|\sum_{i=1}^{k}a_{ni}X_i\right|>4\epsilon\right)\subset\bigcup_{j=1}^{4}\left(\max_{1\le k\le n}\left|\sum_{i=1}^{k}X_{ni}^{(j)}\right|>\epsilon\right).
 \end{eqnarray}
 Thus, in order to establish (\ref{07}), it suffices to prove that
  \begin{eqnarray}\label{14}
I_j:=\sum_{n=1}^{\infty}n^{r-2}\vv\left(\max_{1\le k\le n}\left|\sum_{i=1}^{k}X_{ni}^{(j)}\right|>\epsilon\right)<\infty, j=1,2,3,4.
 \end{eqnarray}

 By the definition of $X_{ni}^{(4)}$, we deduce that $\left(\max_{1\le k\le n}\left|\sum_{i=1}^{k}X_{ni}^{(4)}\right|>\epsilon\right)\subset \left(\max_{1\le i\le n}|a_{ni}X_i|>\epsilon/K\right)$. Since $a_{ni}\approx (i/n)^{\beta}(1/n^p)$, by Lemma \ref{lem1}, we see that
 \begin{eqnarray}\label{15}
\nonumber I_4&\le&\sum_{n=1}^{\infty}n^{r-2}\sum_{i=1}^{n}\vv\left(|a_{ni}X_i|>\epsilon/K\right)\le \sum_{n=1}^{\infty}n^{r-2}\sum_{i=1}^{n}\vv\left(|X|>\frac{\epsilon}{K}n^{p+\beta}i^{-\beta}\right)\\
\nonumber&\approx& \int_{1}^{\infty}x^{r-2}\int_{1}^{x}\vv\left(|X|>\frac{\epsilon}{CK}x^{p+\beta}y^{-\beta}\right)\dif y\dif x\mbox{ ( Letting $u=x^{p+\beta}y^{-\beta}$, $v=y$)}\\
\nonumber&=& \frac{1}{p+\beta}\int_{1}^{\infty}\dif u\int_{1}^{u^{1/p}}u^{(r-1)/(p+\beta)-1}v^{\beta(r-1)/(p+\beta)}\vv\left(|X|>\frac{\epsilon}{CK}u\right)\dif v\\
&\approx&\begin{cases} C\int_{1}^{\infty}u^{\frac{r-1}{p+\beta}-1}\vv\left(|X|>\frac{\epsilon}{CK}u\right)\dif u\ll C_{\vv}\left(|X|^{(r-1)/(p+\beta)}\right),& \text{
for $-p<\beta<-p/r$;}\\
C\int_{1}^{\infty}u^{r/p-1}\ln(u)\vv\left(|X|>\frac{\epsilon}{CK}u\right)\dif u\ll C_{\vv}\left(|X|^{r/p}\ln(1+|X|)\right), & \text{ for $\beta=-p/r$;}\\
C\int_{1}^{\infty}u^{r/p-1}\vv\left(|X|>\frac{\epsilon}{CK}u\right)\dif u\ll C_{\vv}\left(|X|^{r/p}\right) , & \text{ for $\beta>-p/r$;}
\end{cases}
 \end{eqnarray}
Hence, by (\ref{06}), $I_4<\infty$. By the definition of $X_{ni}^{(2)}$, we see that $X_{ni}^{(2)}>0$. We have
 \begin{eqnarray}\label{16}
\nonumber && \vv\left(\max_{1\le k\le n}\left|\sum_{i=1}^{k}X_{ni}^{(2)}\right|>\epsilon\right)\\
\nonumber && =\vv\left(\sum_{i=1}^{n}X_{ni}^{(2)}>\epsilon\right)\\
 &&\le \vv\left( \mbox {there are at least $K$ indices $i\in [1,n]$ such that $a_{ni}X_i>n^{-\delta}$ }\right)\\
 \nonumber &&\le \sum_{1\le i_1<i_2<\cdots<i_K\le n}\prod_{j=1}^{K}\vv\left(a_{ni_{j}}X_{i_{j}}>n^{-\delta}\right)\le \left(\sum_{j=1}^{n}\vv(a_{nj}X>n^{-\delta})\right)^K.
 \end{eqnarray}
 Since (\ref{06}) implies $\ee|X|^{r/p}<\infty$, by Markov's inequality under sublinear expectations and (\ref{16}), we get
  \begin{eqnarray}\label{17}
\nonumber I_2&\le& \sum_{n=1}^{\infty}n^{r-2}\left(\sum_{j=1}^{n}\vv(a_{nj}X>n^{-\delta})\right)^K\\
\nonumber &\le& \sum_{n=1}^{\infty}n^{r-2}\left(\sum_{j=1}^{n}n^{r\delta p}|a_{nj}|^{r/p}\ee|X|^{r/p}\right)^K\\
 &\approx& \begin{cases}\sum_{n=1}^{\infty}n^{r-2-Kr(p+\beta-\delta/p)} ,& \text{
for $-p<\beta<-p/r$;}\\
\sum_{n=1}^{\infty}n^{r-2-K(r-1-r\delta/p)}\log^Kn, & \text{ for $\beta=-p/r$;}\\
\sum_{n=1}^{\infty}n^{r-2-K(r-1-r\delta/p)}, & \text{ for $\beta>-p/r$.}
\end{cases}
 \end{eqnarray}
Notice that $r>1$, $p+\beta>0$, we could choose $\delta>0$ small enough and sufficiently large integer $K$ such that $r-2-Kr(p+\beta-\delta)/p<-1$ and $r-2-K(r-1-r\delta/p)<-1$. Hence, by (\ref{17}), we obtain $I_2<\infty$. Similarly, we could get $I_3<\infty$. In order to estimate $I_1$, we verify that
\begin{eqnarray}\label{18}
\max_{1\le k\le n}\left|\sum_{i=1}^{k}\ee X_{ni}^{(1)}\right|\rightarrow 0 \mbox{   as $n\rightarrow\infty$.}
\end{eqnarray}
Notice that (\ref{06}) implies $\ee|X|^{r/p}<\infty$ and $\ee|X|^{1/p}<\infty$, when $p>1$, note that $|X_{ni}^{(1)}|\le n^{-\delta}$ and $|X_{ni}^{(1)}|\le |a_{ni}X_i|$, we see that
\begin{eqnarray}\label{19}
\nonumber\max_{1\le k\le n}\left|\sum_{i=1}^{k}\ee X_{ni}^{(1)}\right|&\le& \sum_{i=1}^{n}\ee \left|X_{ni}^{(1)}\right|\\
\nonumber &\le&n^{-\delta(1-1/p)}\sum_{i=1}^{n}\ee \left|a_{ni}X_i\right|^{1/p}\\
&\ll&n^{-\delta(1-1/p)}\sum_{i=1}^{n}n^{-(p+\beta)/p}i^{\beta/p}\approx n^{-\delta(1-1/p)}\rightarrow 0 \mbox{ as $n\rightarrow\infty$. }
\end{eqnarray}
When $\frac12<p\le 1$, notice that $\ee[X]=\ee[-X]=0$, by choosing $\delta$ small enough such that $-\delta(1-r/p)+1-r<0$, we get
\begin{eqnarray}\label{20}
\nonumber \max_{1\le k\le n}\left|\sum_{i=1}^{k}\ee X_{ni}^{(1)}\right|&\le&2\sum_{i=1}^{n}\ee|a_{ni}X_i|I(|a_{ni}X_i|>n^{-\delta})\le 2n^{-\delta(1-r/p)}\sum_{i=1}^{n}\ee|a_{ni}X_i|^{r/p}\\
\nonumber &\ll&  n^{-\delta(1-r/p)}\sum_{i=1}^{n}|a_{ni}|^{r/p}\ll n^{-\delta(1-r/p)}\left(\sum_{i=1}^{n}n^{-r(p+\beta)/p}i^{\beta r/p}\right)\\
 &\approx& \begin{cases}n^{\delta-r(p+\beta-\delta)/p} ,& \text{
for $-p<\beta<-p/r$,}\\
n^{\delta(1-r/p)+1-r}\log n, & \text{ for $\beta=-p/r$,}\\
n^{\delta(1-r/p)+1-r}, & \text{ for $\beta>-p/r$,}
\end{cases}\\
&&\rightarrow 0 \mbox{ as $n\rightarrow \infty$.}
 \end{eqnarray}
 Hence, to prove $I_1<\infty$, it suffices to prove that
 \begin{eqnarray}\label{21}
I_1^{*}:=\sum_{n=1}^{\infty}n^{r-2}\vv\left(\max_{1\le k\le n}\left|\sum_{i=1}^{k}\left(X_{ni}^{(1)}-\ee X_{ni}^{(1)}\right)\right|>\epsilon\right)<\infty.
\end{eqnarray}
Notice that $\{X_{ni}^{(1)},1\le i\le n,n\ge 1\}$ is independent, identically distributed. By Markov's inequality under sublinear expectations, $C_r$ inequality, and Lemma \ref{lem3}, we see that for a suitably large $M$, which will be determined as in the proof of Theorem 2.1 in Guo \cite{Guo2012},
 \begin{eqnarray}\label{22}
\nonumber&&\vv\left(\max_{1\le k\le n}\left|\sum_{i=1}^{k}\left(X_{ni}^{(1)}-\ee X_{ni}^{(1)}\right)\right|>\epsilon\right)\\
&&\ll (\log n)^M\left(\sum_{i=1}^{n}\ee\left|X_{ni}^{(1)}\right|^M+\left(\sum_{i=1}^{n}\ee\left|X_{ni}^{(1)}\right|^2\right)^{M/2}\right).
\end{eqnarray}
With Lemma \ref{lem4} here in place of Lemma 1.6 in Guo \cite{Guo2012}, the rest of the proof is similar to that of Theorem 2.1 in Guo \cite{Guo2012}. It is omitted here.

{\bf Proof.}[ Proof of Theorem \ref{thm2}]  The proof is similar to that of Theorem \ref{thm1} and is omitted.  $\Box$

{\bf Proof.}[ Proof of Corollary \ref{cor1}]. Set $a_{ni}=(A_{n-i}^{\alpha-1}/A_{n}^{\alpha})^p$, $0\le i\le n$, $n\ge 1$. Notice that, for $\alpha>-1$, $A_n^{\alpha}\approx n^{\alpha}/\Gamma(\alpha+1)$. Hence, for $\alpha>0$, we get $a_{ni}\approx (n-i)^{p(\alpha-1)}n^{-p\alpha}$, $0\le i\le n$, $n\ge 1$, $a_{nn}\approx n^{-p\alpha}$. Therefore, setting $\beta=p(\alpha-1)$ in Theorem \ref{thm2}, we can see that (\ref{10}) is equivalent to (\ref{11}).

%\Acknowledgements{This work was Supported by grants from the NSF of China (Grant No.11731012,12031005),   Ten Thousands Talents Plan of Zhejiang Province (Grant No. 2018R52042) and the Fundamental Research Funds for the Central Universities.}


\end{document}